\documentclass{amsart}
\usepackage{amssymb}

\newcommand{\calg}{$C^{*}$-algebra\;}
\newcommand{\calgs}{$C^{*}$-algebras\;}
\newcommand{\ts}{topological space\;}
\newcommand{\cpt}{compact\;}
\newcommand{\cgen}{compactly generated\;}
\newcommand{\creg}{completely regular\;}
\newcommand{\hau}{Hausdorff\;}
\newcommand{\chau}{completely Hausdorff\;}
\newcommand{\lcpt}{locally compact\;}
\newcommand{\cnt}{continuous\;}
\newcommand{\ai}{approximate identity\;}
\newcommand{\fn}{function\;}
\newcommand{\fns}{functions\;}
\newcommand{\homo}{homomorphism\;}
\newcommand{\shomo}{$*$-homomorphism\;}
\newcommand{\homos}{homomorphisms\;}
\newcommand{\repn}{representation\;}
\newcommand{\repns}{representations\;}
\newcommand{\pro}{projective limit\;}
\newcommand{\Apro}{$A=\varprojlim_\alpha A_\alpha$\;}
\newcommand{\Asig}{$A=\varprojlim_n A_n$\;}

\newcommand{\Bpro}{$B=\varprojlim_\alpha B_\alpha$\;}
\newcommand{\proc}{pro-$C^*$-algebra\;}
\newcommand{\procs}{pro-$C^*$-algebras\;}
\newcommand{\sig}{$\sigma$-$C^*$-algebra\;}
\newcommand{\sigs}{$\sigma$-$C^*$-algebras\;}
\newcommand{\pros}{pro-$C^*$-algebras\;}

\newtheorem{defi}{Definition}[section]
\newtheorem{prop}{Proposition}[section]
\newtheorem{theo}{Theorem}[section]
\newtheorem{lemm}{Lemma}[section]
\newtheorem{cor}{Corollary}[section]
\newtheorem{ex}{Example}[section]
\newtheorem{rem}{Remark}[section]

\begin{document}
\title[pro-$C^*$-algebras]{Locally compact pro-$C^*$-algebras}
\author[M. Amini]{Massoud Amini}
\address{Department of Mathematics and Statistics\\ University of Saskatchewan
\\106 Wiggins Road, Saskatoon\\ Saskatchewan, Canada S7N 5E6\\mamini@math.usask.ca}
\thanks{Partially supported by a grant from IPM.}
\maketitle

\begin{abstract}
Let $X$ be a \lcpt non \cpt \hau \ts. Consider the algebras
$C(X)$, $C_b(X)$, $C_0(X)$, and $C_{00}(X)$ of respectively
arbitrary, bounded, vanishing at infinity, and compactly supported
\cnt \fns on $X$. From these, the second and third are \calgs, the
forth is a normed algebra, where as the first is only a
topological algebra (it is indeed a \proc ). The interesting fact
about these algebras is that if one of them is given, the rest can
be obtained using functional analysis tools. For instance, given
the \calg $C_0(X)$, one can get the other three algebras by
$C_{00}(X)=K(C_0(X))$, $C_b(X)=M(C_0(X))$,
$C(X)=\Gamma(K(C_0(X)))$, where the right hand sides are the Pedersen's ideal, the multiplier algebra, and the unbounded multiplier algebra of the Pedersen's ideal of $C_0(X)$, respectively. 
In this article we consider the
possibility of these transitions for general \calgs . The
difficult part is to start with a \proc $A$ and to construct a
\calg $A_0$ such that $A=\Gamma(K(A_0))$. The \procs for which
this is possible are called {\it locally compact} and we have
characterized them using a concept similar to \ai .
\end{abstract}

\section{Multipliers of Pedersen's ideal}

The \calg of \cnt \fns vanishing at infinity has always been a
source of motivations for many abstract aspects of the theory of
\calgs. Developing a measure theory for \calgs, Gert K. Pedersen
found a minimal dense ideal inside each \calg which plays the role
of \fns of \cpt support in the commutative case [Pd66]. We use
this ideal and objects related to it frequently in next section.
Therefore we would like to give more details about this ideal
here. G.K. Pedersen in his PhD Thesis in mid $60$'s investigated a
non commutative measure theory for \calgs [Pd64], [Pd66, I-IV].
Because most of the interesting measures are not finite, he had to
choose one of the following options: either to consider (infinite)
measures on the state space of the \calg or to regard these
measures as extended valued linear functionals (now called {\it
weights}) on the \calg. The second approach seems more effective,
but then he had to make sure that these weights would be finite on
a fairly large subalgebra. The Pedersen ideal is one candidate.
Indeed he got more: He proved that every \calg $A$, has a dense
two sided ideal $K(A)$ which is minimal (indeed minimum) among all
dense hereditary ideals of $A$, where hereditary means that any
positive element of $A$ majorized by an element of $K(A)$ actually
belongs to $K(A)$ (Later it was proved that it is indeed minimal
(minimum) among all dense ideals [LS]). Two classical examples are
$A=C_0(X)$, algebra \cnt \fns vanishing at infinity, and $B=K(H)$,
algebra of \cpt operators, for which the Pedersen ideal is
$K(A)=C_{00}(X)$, \cnt \fns of \cpt support, and $K(B)=F(H)$,
algebra of finite rank operators. In general, one can explicitly
construct the Pedersen ideal as follows. For a \calg $A$, let
\begin{eqnarray}
K_0(A)&=&\{x\in A^+: xy=x, \ \ some \ \ y\in A^+\}, \\
K_+(A)&=&\{x\in A^+: \exists n\geq 1\ \exists x_1,\dots,x_n\in K_0(A)\ \
x\leq\sum_{k=1}^n x_k\ \}.
\end{eqnarray}
Then $K(A)=spanK_+ (A)$ is a minimal (the minimum) dense ideal of
$A$. It is well known that there always exist an \ai of $A$ inside any
given dense ideal [Mur]. In particular there is an \ai of $A$ consisting of
elements of $K(A)$. Indeed in {\it separable} case we can even choose
the \ai inside $K_0(A)$. More precisely, any separable \calg $A$ has
a countable \ai $(e_n)$ which is {\it canonical} in the sense that
$e_n\geq 0$ and $e_n e_{n+1}=e_n$, for each $n$ [CF]. In
particular $e_n\in K_0(A)$ by definition.

We said that He got more than what is really needed to support
weights on a \calg , but indeed, like almost any other thing in
Mathematics (in life!?), he had to pay the price. The weight
theory extended based on the Pedersen ideal (called $C^*$-{\it
integrals}) only works for weights which are {\it unitarily
bounded}, namely those linear functionals $f:A^+\to [0,\infty]$
for which $$sup\{f(u^*au): u\in\tilde A\ \text{unitary}\}<\infty\
\ (a\in A).$$ This is a strong condition not satisfied for the
canonical weight (obtained by evaluation at the identity) on the
group algebra $C^*(G)$ of a discrete infinite topological group
$G$ [Pd]. He showed that these integrals correspond to unitarily
bounded positive linear functionals on $K(A)$. These could be
decomposed as a (countable) sum of continuous positive linear
functionals, and be represented as integrals on the Pure state
space of $A$. For the above two examples, the $C^*$-integrals
correspond to the set (lattice) of all (positive) Borel measures
on $X$ and positive bounded operators on $H$, respectively, where
in the second example each $S\in B(H)^+$ give rises to the
$C^*$-integral $$f_S(T)=tr(TS)\ \ (T\in F(H)),$$ and all the
$C^*$-integrals on $K(H)$ are of this form [Pd66, I]. Pedersen
also investigated the case where we are interested in the absolute
value of an integral. For any linear functional on $K(A)$, its
absolute value is a convex functional on $K_{+}(A)$. Considering a
convex functional $f$ on $A^+$ he associated a {\it variation
functional} to it, defined by $$var(f)(a)=sup\{f(y^*y): yy^*\leq
x\}\ \ (a\in A^+).$$ Then $var(f)$ is an {\it invariant} (i.e.
tracial and positive) convex functional on $A^+$. He showed that
there is a one-one correspondence between densely defined, lower
semi continuous, invariant convex functionals on $A^+$ and finite
invariant convex functionals on $K_{+}(A)$. In particular, we are
interested in those functionals $f$ on $K(A)$ for which $var(f)$
is finite (these are called of {\it finite variation}). Then we
know that functionals of finite variations are indeed
$C^*$-integrals [Pd66, III].

Let $J$ be a topological (complex) algebra with approximate
identity. By a ({\it double}) {\it multiplier} (or {\it double
centralizer} ) of $J$ we mean a pair $(S,T)$ of functions from $J$
to $J$ such that $$xS(y)=T(x)y,$$ for all $x,y\in J$. One can then
show that $S$ and $T$ have to be linear and, respectively, a {\it
left} and {\it right multiplier} [LT, 2.10], i.e. $$S(xy)=S(x)y,\
\ T(xy)=xT(y)\ \ \ (x,y\in J).$$ We denote the set of all (left,
right) multipliers of $J$ by $\Gamma(J)$ ($\Gamma_\ell (J)$,
$\Gamma_r (J)$, respectively). Then $\Gamma(J)$ is a vector space
under the natural operations, and an algebra under the
multiplication $(S,T)(U,V)=(SU,VT)$. If $J$ is a $*$-algebra, then
so is $\Gamma(J)$, under the involution $(S,T)^* =(S^* ,T^* )$,
where $S^* (x)=S(x^* )^*$, for each $x\in J$. It is obvious that
$\Gamma(J)\subseteq \Gamma_\ell (J)\times \Gamma_r (J)$. If $J$ is
complete, all these three spaces could be identified with
subspaces of $J^{**}$, and $\Gamma(J)=M(J)=\Gamma_\ell
(J)\cap\Gamma_r (J)$. If $J$ is a normed algebra with a
contractive approximate identity, and $(S,T)\in \Gamma(J)$, then
$S$ is bounded if and only if $T$ is bounded and $\| S\|=\| T\|$.
If this is the case, we say that $(S,T)$ is {\it bounded} and put
$\| (S,T)\|=\| S\|$. We denote the set of all bounded (left,
right) multipliers of $J$ by $M(J)$ ($M_\ell (J)$, $M_r (J)$,
respectively). If $J$ is complete, then $\Gamma(J)=M(J)$.

Let $A_0$ be a \calg and $A_{00}=K(A_0 )$ be its Pedersen ideal.
We are particularly interested in $\Gamma (A_{00})$. One reason is
that this algebra contains some unbounded elements which are
important in functional analytic applications. $\Gamma (A_{00})$
has been extensively studied by A.J. Lazar and D.C. Taylor [LT]
and N.C. Phillips [Ph88b]. We sketch some of the main results from
these references.

Every bounded multiplier on $  A_{00}$ extends
uniquely to a multiplier on $A_0$. $A_0$ could be naturally
embedded in $\Gamma (A_{00})$ by identifying
$a\in A_0$ with $(S_a , T_a)\in \Gamma (A_{00})$ defined by
$$S_a (b)=ab, \ T_a (b)=ba.$$

It was first pointed out by N.C. Phillips that $\Gamma (A_{00})$ could
be represented as a projective limit of unital \calgs. Indeed
$$\Gamma (A_{00})=\varprojlim_{I} M(I),$$
where $I$ runs over all closed two sided ideals of $A_0$ singly
generated by elements of (an approximate identity contained in)
$A_{00}^+$ [Ph88b].

There are dif and only iferent topologies on $\Gamma(A_{00})$.
Lazar and Taylor introduced the $\kappa$-{\it topology} on
$\Gamma(A_{00})$. this is the topology induced by seminorms
$$(S,T)\mapsto \|S(a)\|\ \ \text{and} \ \ (S,T)\mapsto \|T(a)\|,$$
where $a$ runs over $A_{00}$. A sequence $(x_i)\subseteq
\Gamma(A_{00})$ converges to $x\in  \Gamma(A_{00})$ if and only if
$\|x_i a-xa\|\to 0$ and $\|ax_i -ax\|\to 0$, as $i\to\infty$, for
each $a\in A_{00}$. Also by Phillips' result we have the {\it
projective topology} induced by the norm topologies of \calgs
$M(I)$ on $\Gamma(A_{00})$. It follows from [Ph88b, thm4] that the
projective topology is stronger than the $\kappa$-topology.

It was already noticed by Lazar and Taylor that $\Gamma (A_{00})$
could be represented as an algebra of
densely defined linear operators on a Hilbert space. We describe this
in some more details. Let $A_0$ be a \calg which is non degenerately
represented on a Hilbert space $H$ (i.e. the closed linear span
of $A_0 H=\{a\zeta :a\in A_0 ,\ \zeta\in H\}$ is equal to $H$). Let $H'
=span(A_{00}H)$. This is clearly a
dense subspace of $H$. We introduce an
algebra of operators on $H'$ as follows
$$B(A_{00},H' )=\{T\in L(H{'} ): xT\ \text{and}\ Tx\in B(H' )\cap A_{00}\
\ (x\in A_{00}) \}$$
where we have identified elements of $B(H' )$ with their unique
extension to $H$ (This is from [LT], but they don't use this
notation).
This is obviously an algebra under the usual
operations. Take $T\in L(H' )$ and let $T^*$ be the (possibly
unbounded) adjoint of this densely defined operator restricted to
$H'$. Let $x\in A_{00}, \zeta\in H$, and $\eta\in H' $, then $y=x^*
T\in A_{00}$ and $<T\eta ,x\zeta>=<\eta , y^* \zeta >$, so $T^*$ is
defined on the whole $H'$ and $T^*(x\zeta)=(x^* T)\zeta$. Under
this involution $B(A_{00},H' )$ becomes a $*$-algebra. Then one can
show that $\Gamma(A_{00})\cong B(A_{00},H' )$
as $*$-algebras [LT].

Indeed to each $u\in\Gamma(A_{00})$ there corresponds $\hat u \in L(H')$
defined by $\hat u (x\zeta)=(ux)\zeta$. Since there is an
approximate identity for $A_0$ in $A_{00}$, $\hat u$ is well
defined. Now for each $x\in A_{00}$, $x\hat u +\hat u x$ is the
restriction of $xu+ux$ to $H'$, and so $\hat u\in B(A_{00},H'
)$. On the other hand, for each $T\in B(A_{00},H')$, define $S(x)$
and $T(x)$ for $x\in A_{00}$ to be extensions of $Tx$ and $xT$ to $H$
and observe that $\hat u =T$, for $u=(S,T)\in \Gamma(A_{00})$. It is
now easy to check that $u\mapsto \hat u$ is a $*$-isomorphism.

Let's make an immediate observation: If $span(A_{00}H)=H$, then
$\Gamma(A_{00})=M(A_0)$. Indeed, in this case,
for each $T\in\Gamma(A_{00})$, $T^*$ is
defined everywhere, and so bounded. Therefore $T$ has to be bounded also.

\section{Locally \cpt \sigs}

We were aiming at the problem of transitions between the
categories of \calgs and \pros . Before dealing with this problem,
however, we should overcome a technical difficulty. Starting with
a \calg $A_0$, one can simply get the other three algebras by
putting $A_{00}=K(A_0), A_b=M(A_0)$, and $A=\Gamma(A_{00})$, but
starting with a \proc $A$ it is not clear how we can appropriately
associate a \calg $A_0$ to it (specially if we want them to be
related by the relation $A=\Gamma(K(A_0))$).
 Indeed, even in commutative case, this is not possible in
general. The commutative unital \sigs are exactly the algebras
$C(X)$ of all \cnt \fns on a countably \cgen \hau \ts $X$, with
the compact-open topology [Ph88a, 5.7]. But $X$ need not to be
\lcpt and so $C_0(X)$ is not necessarily a \calg. However, we show
that it is possible to distinguish a subcategory of \sigs for
which one can make the desired transition. These are naturally
called {\it locally compact \sigs}. This subcategory would contain
all unital \calgs.

Let's first consider the commutative case. For a (\lcpt)
topological space $X$, the relation between compactness type
conditions on $X$ and existence of special types of approximate
identities in $C_0(X)$ has been extensively studied. Here we quote
some of these results to motivate our approach. A countable
approximate identity $(e_n )$ in a Banach $*$-algebra is called
{\it well-behaved} if, for each $n$ and each strictly increasing
subsequence of indices $(n_i)$, there is $N\geq 0$ such that $e_n
e_{n_k}=e_n e_{n_l}\ \ (k,l\geq N)$ [Ty72], [CD]. The compactness
type conditions on a \lcpt \ts $X$ are related to the existence of
certain approximate identities on $C_0 (X)$. H.S. Collins and J.R.
Dorroh showed that $X$ is $\sigma$-\cpt if and only if $C_0 (X)$
has a countable canonical (in the sense of Definition 2.5)
approximate identity [CD]. Collins and Fontenot showed that $X$ is
paracompact if and only if $C_0 (X)$ has a strictly totally
bounded (canonical) approximate identity [CF], and conjectured
that this is also equivalent to the condition that $C_0 (X)$ has a
well-behaved approximate identity (this was proved by R.F. Wheeler
[Whe]). They also showed that if $X$ is pseudo-\cpt and $C_0 (X)$
has a well-behaved approximate identity, then $X$ is compact.
Later, Fontenot and Wheeler showed that $X$ is paracompact if and
only if $C_0 (X)$ has a weakly compact approximate identity [FW].

Now let $A=\varprojlim_{n} A_n$
be a unital (with unit 1) \sig and $\pi_n:
A\to A_n$ be the corresponding \homos and $A_b$ the
(unital) \calg of bounded elements of $A$, that is
$$A_b=\{a\in A: sup_n \|\pi_n (a)\|_n <\infty\}.$$
\begin{defi} An element $x\in A$ is called strongly bounded if
$AxA\subseteq A_b$. We denote the set of all such
elements by $A_{sb}$.
\end{defi}

Then $A_{sb}\subseteq A_b$ is clearly a two sided ideal of $A$.
\label{mb}
\begin{defi}  An element $x\in A^+$ is called
 multiplicatively bounded (m.b.) if $Ax\subseteq A_b$. A sequence in
$A^+$ is called multiplicatively bounded if all of its elements are
multiplicatively bounded.
\end{defi}

Note that in the above definition,
since $x$ is positive (and so self adjoint), then
the above condition would imply $xA\subseteq A_b$.
In particular $x^2\in A_{sb}$.
\label{support}
\begin{defi}
Let $A$ be a \sig . A {\it support algebra} of $A$ is a two sided
ideal $A_{00}$ of $A$ of the form $$A_{00}=\bigcup_n Ae_n ^2 A ,$$
where $e_n$'s are some given multiplicatively bounded elements of
$A^+$ (compare with [Lin]).
\end{defi}

Each support algebra of $A$ is clearly a two sided ideal of $A$.
Also it is always contained in $A_{sb}$. This ideal then induces a
topology on $A$ as follows.
\begin{defi} The $A_{00}$-topology on $A$ is the weakest topology
such that the maps from $A$ to $A_b$ of the form $$a\mapsto xa\ \
\text{and} \ \ a\mapsto xa,$$ are \cnt for each $x\in A_{00}$,
where $A_b$ has its norm topology. A sequence $(a_i)\subseteq A$
converges to $a\in A$ in the $A_{00}$-topology if and only if
$\|xa_i -xa\|_\infty\to 0$ and $\|a_i x-ax\|_\infty\to 0$, as
$i\to\infty$, for each $x\in A_{00}$.
\end{defi}
\label{aid}
\begin{prop} Let $(e_n)$ be a m.b. sequence in $A$ and let $A_{00}$
be the corresponding support algebra of $A$. Consider the
following conditions:

$(1)$ $(e_n)$ is an \ai of $A$ in the $A_{00}$-topology,

$(2)$ $e_n\to 1$ in the $A_{00}$-topology,

$(3)$ $(e_n)$ is an \ai of $A_{00}$ in the norm topology.

Then $(1)\Rightarrow  (2) \Leftrightarrow (3)$.

Moreover if
$$\forall n\forall a\in A\ \ \|\pi_n(a)\|_n\geq\|e_na\|_\infty ,$$
then every \ai in the projective topology is an \ai in the
$A_{00}$-topology.
\end{prop}

{\bf Proof} The fact that $(1)$ implies $(2)$ is immediate (as $1\in
A$). The equivalence of $(2)$ and $(3)$ is just the definition of
the $A_{00}$-topology. The last statement follows from the fact that
condition mentioned above says that the projective topology is stronger than
the $A_{00}$-topology. Indeed, if $(x_i)$ be a sequence in $A$ and $x_i\to
x\in A$ in the projective topology, then for each $n$ and $a,b\in A$
we have
$$\|ae_n ^2 b(x_i -x)\|_\infty \leq\|ae_n\|_\infty\|\pi_n(b)\|_n \|\pi_n (x_i
-x)\|_n \to 0,$$
as $i\to \infty$. Now each element of $A_{00}$ is a finite sum of
elements of the form $ae_n b$, and we conclude that $x_i\to x$ in
the $A_{00}$-topology.\qed
\label{canonical}
\begin{defi} Let $(e_n)\subseteq A^+$ be a m.b. sequence in
$A$ and let $A_{00}$ be the corresponding support algebra. Then
$(e_n)$ is called canonical if $e_m e_n=e_n$, for each $n<m$.
\end{defi}

\label{lcpt}
\begin{defi} A unital \sig $A$ is called {\it locally compact} if it has a
presentation
$A=\varprojlim_{n} A_n$
as an inverse limit of a countable family of
\calgs such that the corresponding morphisms $\pi_n :A\to A_n$
are surjective, and
there is a multiplicatively bounded, canonical
sequence $(e_n)_{n\in \Bbb N}$ with $0\leq
e_n \leq 1\ \ (n\in \Bbb N)$, such that if $A_{00}=\bigcup_n
Ae_n ^2 A $ is the corresponding support algebra of $A$, then
$e_n\to 1$ in the $A_{00}$-topology, $A$ is complete in the
$A_{00}$-topology, and it
satisfies the following "compatibility condition":
$$\forall n\ \forall a\in A\ \ sup\{\|e_n ba\|_\infty :\ b\in A,\ \|e_n
b\|_\infty \leq 1\}\geq\|\pi_n (a)\|_n \geq \|e_na\|_\infty.$$
\end{defi}

Note that the above condition implies its right hand side version
(i.e. with $e_n$ multiplied from right). This follows from the fact
that $\pi_n$ is involutive and $A$ is a $*$-algebra.


In order to get a better idea about these definitions, let's look at
some classical examples.
\label{worse}
\begin{ex} Let $A=C(\Bbb R)$ be the \sig of all \cnt \fns on the real
line. For each $n\geq 1$, let $A_n=C[-n,n]$ be the unital \calg of
all \cnt \fns on the \cpt interval $[-n,n]$. For $m\geq n$ we have
the morphism $\pi_{nm}:A_m \to A_n$ which sends a \cnt \fn on
$[-m,m]$ to its restriction on $[-n,n]$. Then $(A_n,\pi_{nm})$ is
an inverse system of \calgs with $$C(\Bbb R)=\varprojlim_n
C[-n,n],$$ where the morphism $\pi_n :C(\Bbb R)\to C[-n,n]$ is
also defined by restriction. Then, for $f\in C(\Bbb R)$, $\|\pi_n
(f)\|=\|f|_{[-n,n]}\|$, so $\|f\|_\infty$ is nothing but the usual
sup-norm. In particular, the set of bounded elements of $A$ is
exactly the (unital) \calg $A_b=C_b(\Bbb R)$ of bounded \cnt \fns
. Also it is easy to see that the strongly bounded elements are
exactly those of \cpt support, namely $A_{sb}=C_{00}(\Bbb R)$.
Also a sequence $(f_n)\subseteq C(\Bbb R)$ converges to $f\in
C(\Bbb R)$ in the projective topology if and only if it uniformly
converges to $f$ on each interval $[-n,n]$. Since each \cpt subset
of $\Bbb R$ is already contained in one of these intervals, the
projective topology is nothing but the topology of uniform
convergence on \cpt sets (or the so called {\it compact-open}
topology).

Now consider the sequence $(e_n)$ in $C(\Bbb R)$, where $0\leq
e_n\leq 1$ is the {\it bump \fn} which is $1$ on $[1-n,n-1]$ and
zero off $[-n,n]$. Let's observe that this is an \ai of $C(\Bbb
R)$ in the projective topology. Given $m\geq 1$ and $f\in C(\Bbb
R)$, we need to show that $e_n f$ converges to $f$ uniformly on
$[-m,m]$. But this becomes obvious when we note that $e_n=1$ on
$[-m,m]$ for $n\geq m+1$. Now as $e_n$'s have \cpt support, this
sequence is multiplicatively bounded. For each $n$, $Ae_n$
consists of those elements of $A$ whose support is inside the
support of $e_n$, so the corresponding support algebra $A_{00}$ of
$A$ is nothing but $C_{00}(\Bbb R)$. Hence the $A_{00}$-topology
is also the compact-open topology and the two topologies coincide.
Also for $n>m$, $e_n$ is $1$ on the support of $e_m$, so $e_n
e_m=e_m$, i.e. we have a canonical \ai . Also, given $f\in C(\Bbb
R)$ and $n\geq 1$, $e_n$ is dominated by the characteristic \fn of
the interval $[-n,n]$, so in particular $\|fe_n\|_\infty\leq
\|f|_{[-n,n]}\|_\infty$. The other inequality in the compatibility
condition holds for $e_n$ replaced by $e_{n+1}$ (see the remark
after Theorem 3.1).
\end{ex}

\begin{ex} Let $A_n=\ell_n ^\infty$ be
the direct sum of $n$ copies of $\Bbb C$ and for $m\geq n$, let
$\pi_{nm}: A_m\to A_n$ be the projection onto the first $n$
components. This is an inverse system of \calgs with the inverse
limit being the algebra $A=c$ of all sequences. Then $\pi_n: A\to
A_n$ is simply the projection onto the first $n$ components. Also
$A_b=\ell ^\infty$ and $A_{sb}$ is the set of all sequences with
finitely many non zero components. The projective topology on $A$
is just the topology of pointwise convergence, a sequence of
sequences converges to a sequence in projective topology if and
only if for each $n$ the sequence of the $n$-th components
converges to the $n$-th component in $\Bbb C$. Here
$e_n=(1,\dots,1,0,0,\dots)$ with the first $n$ components $1$ and
the rest $0$, forms a multiplicatively bounded canonical \ai of
$A$ in projective topology which satisfies the compatibility
condition. The corresponding support algebra is $A_{00}=A_{sb}$
and the $A_{00}$ topology is the same as the projective topology.
\end{ex}
\begin{ex} Let $A=B(\ell ^2)$ be the \calg of all bounded operators
on the separable Hilbert space $\ell ^2$. This is a unital \sig in a
trivial way (it is the inverse limit of the inverse system consisting
of only one unital \calg , namely $B(\ell ^2)$ itself with the
connecting homomorphism being identity!). In this case every element is
(strongly) bounded, so $A_b=A_{sb}=A$. The projective topology
coincides with the norm topology, so an \ai in projective topology
is just a sequence
which converges to $I$ in norm. Since all elements are bounded in this
case, every \ai is multiplicatively bounded. The trivial \ai
consisting only of the identity operator $I$ is canonical and
satisfies the compatibility condition.  This choice of \ai gives
$A_{00}=A$ and the two topologies coincide.

There is, however, a more interesting choice of \ai . Take the
presentation $B(\ell ^2)=\varprojlim_n A_n$, where $A_n=B(\ell^2)$
and $\pi_{mn}:A_m\to A_n$ is the identity,
for each $m,n\geq 1$. Then again $A_b=A_{sb}=A$ and the projective
topology is the norm topology.
For each $n$, put $e_n= \begin{pmatrix} I_n&0\\ 0&0\end{pmatrix}\in
 B(\ell^2)$. then
$0\leq e_n\leq I$ and $(e_n)$ is a canonical multiplicatively
bounded sequence in $B(\ell^2)$. The corresponding support algebra
$A_{00}=\bigcup_n  B(\ell^2)e_n B(\ell^2)=F(\ell^2)$ is the
algebra of finite rank operators. Hence the $A_{00}$ topology is
the strong${}^*$-topology on $B(\ell^2)$. Now note that $(e_n)$ is
an \ai of $B(\ell^2)$ only in this weaker topology. To see that
the compatibility condition holds in this case, one only needs to
observe that $\pi_n=id$ and, given an infinite matrix $a\in
B(\ell^2)$, $e_m a$ is obtained from $a$ by fixing the first $m$
rows and making the rest zero.
\end{ex}

\begin{ex} Let $\Bbb M_n=\Bbb M(n,\Bbb C)$ be the unital
\calg of all $n\times
n$ complex matrices, and
$$A_n=\Bbb M_{1,\dots,n}=\sum_{k=1}^n \bigoplus \Bbb M_k ,$$
be the \calg of all block matrices with increasing blocks of size from
$1$ to $n$. Then for each $m\geq n$, the projection $\pi_{nm}: A_m\to A_n$
on the first $n$ block of the top left hand corner is a morphism of \calgs
, and $(A_n,\pi_{nm})$ is an inverse system of \calgs . The
corresponding \sig
$$A=M_{\Bbb N}=\varprojlim_n
 \sum_{k=1}^n \bigoplus \Bbb M_k =\prod_n \Bbb M_n ,$$
is the algebra of all infinite block matrices with an increasing
sequence of blocks of size $1,2,\dots$. Also $\pi_n :A\to A_n$ is
simply the projection onto the first $n$ block. A typical element
of $A$ is of the form $M=(M_1,\dots,M_n,\dots)$, where $M_n\in
\Bbb M_n$. Then $\pi_n(M)=(M_1,\dots,M_n)$ has norm $sup_{1\leq
i\leq n} \|M_i\|$. Therefore $A_b$ consists of those block
matrices $M$ for which $sup_i \|M_i\|<\infty$. Also $A_{sb}$
consists of the block matrices with only finitely many non zero
blocks. A sequence of infinite block matrices converges in the
projective topology to an infinite block matrix if and only if for
each $n$, the $n$-th block of the elements of the sequence
converges to the $n$-th block of that element in norm. If $I_n$
denote the $n\times n$ identity matrix, then the sequence $(e_n)$,
where $e_n=(I_1,I_2,\dots,I_n,0,0,\dots)$, forms a
multiplicatively bounded, canonical \ai of $A$ in the projective
topology which satisfies the compatibility condition. Here
$A_{00}=A_{sb}$ and the two topologies coincide.
\end{ex}

Also one expects that commutative \lcpt unital \sigs should be of the form
$C(X)$ for a $\sigma$-compact , locally compact, (completely) Hausdorff topological 
space $X$. We show that this is true. Let's first see what is the general form of a
commutative unital (not necessarily \lcpt )\,\sigs . We know that each
commutative unital \sigs is of the form $C(X)$ for some
\hau $\kappa_\omega$-space $X$. (recall that $X$ is called a
$\kappa_\omega$-space , or countably \cgen , if it is the union of an
increasing sequence of \cpt subsets which determine the topology).
Indeed the functor
$X\mapsto C(X)$ is a contravariant equivalence. Moreover a \hau
$\kappa_\omega$-space is automatically \creg ( and so \chau) [Ph88a,
prop 5.7]. Also each $\kappa_\omega$-space is $\sigma$-\cpt (by
definition), but $\sigma$-\cpt $\kappa$-spaces are not
$\kappa_\omega$-space necessarily [Ph88a, 5.8].
\label{Xlcpt}
\begin{lemm} If $X$ is a $\sigma$-\cpt \hau \ts and
the $\sigma$-$C^{*}$-algebra $C(X)$ has an approximate identity in the pointwise
convergence topology
consisting of functions of \cpt support, then $X$ is \lcpt.
\end{lemm}
{\bf Proof} Let $(f_n)\subseteq C_{00}(X)$ be an \ai of $C(X)$
in the pointwise convergence topology. Let $F_n$ be the support of $f_n$
and let $V_n=int(F_n)$ be the interior of $F_n$. Since each $V_n$ is an
open set with \cpt closure, it is enough to show that
$\bigcup_n V_n=X$. Let $x\in X$ and let $f\in C(X)$ be the constant
\fn $1$. Then $f_n(x)=
(f_n f)(x)\to f(x)=1$ as
$n\to \infty$. Hence, there is $n$ with
$f_n(x)\geq {\frac  1 2}$, i.e. $x\in V_n$.
This proves the claim and finishes the proof.\qed

\label{cpt}
\begin{lemm} Let $X$ be a \creg \ts and let $K$ be a
\cpt subset of $X$.
Let $f\in C(X)$ have the property that $fg$ is bounded for each
$g\in C(X)$ and $\|fg\|_\infty \leq \|g|_{K}\|_\infty$. Then
$f$ has compact support.
\end{lemm}
{\bf Proof} Let $F=cl(V)$ be the support of $f$,
where $V=\{x\in X: f(x)\neq 0\}$.
Since $X$ is \hau , every closed subset of a \cpt set is \cpt,
therefore we only need to show that $V\subseteq K$. Assume that, on the
contrary, the open set $V\backslash K$ is nonempty. Choose $x\in
V\backslash K$, then there is a \fn $g\in C(X)$ such that $g(x)=1$ and
$g=0$ off $V\backslash K$ ($X$ is \creg). In
particular, $g=0$ on $K$ and so $\|fg\|_\infty\leq\|g|_
{K}\|_\infty=0$. Hence $fg=0$
everywhere. But $fg(x)=f(x)g(x)=f(x)\neq 0$, as $x\in V$, which is a
contradiction. \qed

\label{cor}
\begin{theo} Let $X$ be a \hau $\kappa_\omega$-space
. Then the commutative \sig $C(X)$ is \lcpt if and only if $X$ is
\lcpt.
\end{theo}
{\bf Proof} Let $X$ be \lcpt . We already know that $C(X)=\varprojlim_n
C(K_n)$ for an increasing (ordered by
inclusion) sequence $(K_n)$ of \cpt
subsets of $X$ satisfying $X=\bigcup_n K_n$
with $\pi_n (f)=f|_{K_n}\ \ (f\in
C(X))$ [Ph88a] . Since $X$ is \lcpt we may assume that there are open
sets $V_n$ with \cpt closure such that $K_n\subseteq
V_{n+1}\subseteq K_{n+1}$, for each $n$. (Take an open covering of
$K_n$ by open sets with \cpt closure and choose a finite subcover. Let
$K_n ^{'}$ be the closure of the union of the elements of this
subcover, then $K_n ^{'}$'s are \cpt neighborhoods whose union is $X$.
Next to make this sequence increasing, one can replace $K_n ^{'}$ with $K_n
^{''}=\cup_1 ^n K_i ^{'}$). By
the Urysohn lemma, for each $n\in\Bbb N$, there is an element
$e_n \in C_{00}(X)$ such that $0\leq e_n \leq 1$ and
$e_n =1$ on $K_n$ and $e_n=0$ off $V_{n+1}$.
This is an approximate identity for $C_0(X)$ which satisfies
all of the conditions in Definition~\ref{lcpt}, so
$C(X)$ is a \lcpt \sig.

Conversely let $C(X)=\varprojlim_n C(K_n)$ be \lcpt as a \sig,
where $\{K_n\}$ is an increasing sequence of \cpt subsets of $X$
which determine the topology. Then $C(X)$ has an approximate
identity in compact-open topology (and so in pointwise convergence
topology) which satisfies the compatibility condition of
Definition~\ref{lcpt}. Now $X$ is a \hau $\kappa_\omega$-space and
so it is also \creg and $\sigma$-\cpt. Now by Lemma 2.2, each
member of this \ai has \cpt support, and so $X$ is \lcpt by Lemma
2.1 .\qed

\vspace{.3 cm}
Now let $A$ be a \lcpt \sig , let $(e_n)$
be an \ai of $A$ satisfying the conditions of the above definition, and let
$A_{00}$ be the corresponding support algebra of $A$. Let
$A_0$ be the norm closure of $A_{00}$ in $A_b$. Then $A_0$ is a
(possibly non unital) $C^*$-subalgebra of $A_b$.
We want to show that indeed
$A_{00}=K(A_0)$. First we need some lemmas.

\label{same}
\begin{lemm} With above notations, for each $n$
$$Ae_n ^2 A=A_0 e_n ^2A_0=A_{00} e_n ^2 A_{00}.$$
\end{lemm}
{\bf Proof} For each $n$ we have $$Ae_n ^2
A=Ae_{n+1} e_n ^2 e_{n+1} A\subseteq AA_{00}e_n ^2 A_{00}A\subseteq
A_{00} e_n ^2 A_{00}\subseteq A_0 e_n ^2
A_0\subseteq Ae_n ^2 A,$$
and so all three algebras are equal.\qed
\label{ai}
\begin{lemm}
Let $A, A_{00}, A_0$ and $(e_n)$ be as above.
Then $(e_n)$ is an approximate
identity of $A_0$ in the norm topology.
\end{lemm}
{\bf Proof} We know that $e_n\to 1$ in the $A_{00}$-topology.
Hence, by Proposition 2.1 it forms an \ai for $A_{00}$. Now recall
that $e_n\in A_b$ and $0\leq e_n \leq 1$. By this and the fact
that $A_{00}\subseteq A_0$ is norm dense, the result follows
immediately. Indeed, if $a\in A_0$ and $\epsilon>0$ are given,
then there is $b\in A_{00}$ such that $\|a-b\|_{\infty}
<\epsilon$. Now $b\in A_{00}$ and so by above Lemma, $$\exists n_0
\forall n\geq n_0 \ \ \|be_n -b\|_\infty<\epsilon .$$ Also $\|be_n
-ae_n \|_\infty \leq \|\pi_n(b-a)\|_n \leq \|a-b\|_{\infty}
<\epsilon$. Combining these inequalities we get $\|ae_n
-a\|_\infty < 3\epsilon$, for each $n\geq n_0$.
Hence $(e_n)$ is a right
approximate identity for $A_0$. Similarly one can show that it is also
a left approximate identity.\qed

\vspace{.3 cm}
Next we show that $A_{00}$ is indeed the Pedersen ideal of $A_0$.

\label{K(Ao)}
\begin{lemm} With above notations, $K(A_0)=A_{00}$.
\end{lemm}
{\bf Proof} Since $A_{00}$ is a dense ideal of $A_0$, we have
$K(A_0)\subseteq A_{00}$. Now given $n$, by the fact that $(e_n)$
is canonical, $e_n e_{n+1}=e_n$ and hence $e_n\in
K_{0}(A_0)\subseteq K(A_0)$ . But then, as $K(A_0)\subseteq A_0$
is a two sided ideal, $A_0 e_n\cup e_n A_0\subseteq K(A_0)$ .
Hence, by Lemma 2.3, $A_{00}=\bigcup_n Ae_n ^2 A =\bigcup_n A_0
e_n ^2 A_0 \subseteq K(A_0)$.\qed

\vspace{.3 cm}
Lets check the above examples to see what is $A_0$ in each
example.

\begin{ex}
Let $A=C(\Bbb R)$, then  $A_b=C_b(\Bbb R)$ and for $0\leq e_n\leq
1$ being the {\it bump \fn} which is $1$ on $[1-n,n-1]$ and zero
off $[-n,n]$, we get $A_{00}=C_{00}(\Bbb R)$. Hence
$A_0=C_{0}(\Bbb R)$.
\end{ex}
\begin{ex} Consider the \sig $A=c$ of all sequences of complex numbers,
then $A_b=\ell ^\infty$, and for $e_n=(1,\dots,1,0,0,\dots)$ we get
$A_{00}$ as the ideal of all sequences with
finitely many non zero components. Hence $A_0=c_0$.
\end{ex}
\begin{ex}  Let $A=B(\ell ^2)$ with norm topology. Then $A_b=A$. If
the \ai consists of the identity $I$ only, then $A_{00}=A_{0}=A$. In
the case that we choose the \ai consisting of the elements
$$e_n=\begin{pmatrix} I_n&0\\ 0&0\end{pmatrix} \in B(\ell^2).$$
Then we get $A_{00}=F(\ell^2)$ and so
$A_0=K(\ell^2)$.
\end{ex}
\begin{ex}  For the \sig $A=\prod_n \Bbb
M_n $, we have $A_b=\sum_n ^{\ell^\infty} \bigoplus \Bbb M_n$ and,
if we choose the \ai $e_n=(I_1,I_2,\dots,I_n,0,0,\dots)$, then
$A_{00}=\bigcup_n (\sum_{i=1} ^{n} \bigoplus \Bbb M_i)$ consists
of the block matrices with only finitely many non zero blocks.
Hence $A_0=\sum_n ^{c_0} \bigoplus \Bbb M_n$.
\end{ex}

Note that what is common between these examples is that
$A=\Gamma(K(A_0))$.
Indeed as a typical example of \lcpt \sigs we show that multiplier algebra of
the Pedersen ideal of any \calg is a \lcpt \sig . Before that we
need a trivial lemma which we state it
without proof.

\begin{lemm} Let $X$ be a Banach space and let $X_0\subseteq X$ be a dense
subspace. Let $T_0:X_0\to X_0$ be a bounded linear map. Then $T_0$
extends uniquely to a bounded linear map $T\in B(X)$ and
$$\|T\|=sup\{\|T_0(x)\|: x\in X_0,\ \|x\|\leq 1\}.\qed$$ 
\end{lemm}

Let $A_0$ be a \calg and $A_{00}=K(A_0)$. For each $a\in A_0$, we put
$L_a=cl(A_0 a), \,  R_a =cl(aA_0)$, and $I_a=cl(L_a R_{a^*})$, where the
closures are taken in the norm topology of $A_0$. Consider all
linear maps $S:L_a\to L_a$ and $T:R_{a^*}\to R_{a^*}$ and put
$$M_a=\{(S,T): yS(x)=T(y)x for each \ x\in L_a, y\in R_{a^*}\}.$$
Then $S, T$ are automatically bounded and $M_a$ is a \calg. Indeed
$M_a\cong M(I_a)$, for each $a\in A_0^+$.

\label{yes}
\begin{theo} For each $\sigma$-unital \calg $A_0$, the
\sig $A=\Gamma(K(A_0))$ is locally compact.
\end{theo}
{\bf Proof} Recall that for each (countable contractive) canonical
approximate identity $(e_n)$ of $A_0$ contained in $K(A_0)^+$ (we
can always choose $(e_n)$ to be canonical [CF]) we have
$\Gamma(K(A_0))=\varprojlim_{n} M_{e_n}$ such that the
corresponding projective limit $\pi_n :A\to A_n=M_{e_n}$ sends
$(S,T)\in  \Gamma(K(A_0))$ to $(S|_{L_{e_n}},T|_{R_{e_n}})$. Also
$A_b=M(A_0)$ [Ph88b]. Let $A_{00}$ be the corresponding support
algebra of $A=\Gamma(K(A_0))$; that is $A_{00}=\bigcup_n Ae_n ^2
A$. Then we claim that $A_{00}=K(A_0)$.

Indeed, $K(A_0)$ is an ideal in $A$ [LT], so $A_{00}=\bigcup_n
Ae_n ^2 A \subseteq AK(A_0)A\subseteq K(A_0)$. On the other hand,
$A_{00}$ contains the \ai $(e_n ^2)$ of $A_0$, so it is dense in
$A_0$ (in the norm topology). Being a dense two sided ideal, then
$A_{00}\supseteq K(A_0)$, and the claim is proved. In particular
the $A_{00}$-topology on $A$ is just the $\kappa$-topology. The
fact that $(e_n)$ is an \ai of $A_0$ then implies that $e_n\to
1\in A$ in the $\kappa$-topology. We know that $A$ is complete in
the $\kappa$-topology [LT, 3.8]. We can always choose $(e_n)$ to
be canonical [CF]. The fact that $(e_n)$ is multiplicatively
bounded follows immediately, as for each $n$, $e_n A\subseteq
K(A_0)A\subseteq K(A_0)\subseteq A_b$.

Finally for the
compatibility conditions observe that, for each $n$ and each
$x=(S,T)\in \Gamma(K(A_0))$, we have
$$\|\pi_n((S,T))\|_n=\|(S|_{L_{e_n}},T|_{R_{e_n}})\|_n=\|S|_{L_{e_n}}\|.$$
Now recall that $Ae_n$ is a dense linear subspace of
$L_{e_n}$ and so by the above Lemma, we have
\begin{eqnarray}
\|S|_{L_{e_n}}\| &=sup\{\|S(be_n)\|_\infty :\ b\in A,\ \|be_n\|_\infty\leq 1\}
\\ &=sup\{\|xbe_n\|_\infty :\ b\in A,\ \|be_n\|_\infty\leq 1\}.
\end{eqnarray}

But $e_{n+1}e_n=e_n$ and $\|e_{n+1}e_n\|_\infty\leq 1$, and so the
last supremum is clearly not less than $\|xe_{n+1}e_n\|_\infty =\|xe_n\|_\infty$,
which proves one inequality of the compatibility condition. For
the other, the above calculation shows that indeed equality holds.
\qed

\vspace{.3 cm}
In the next section we would show that the converse is also true.


\vspace{.5 cm}
\section{A covariant functor between two categories}

In this section we show that there is a covariant functor from the
category of $\sigma$-unital \calgs to the category of \lcpt \sigs
. We start with a non unital $\sigma$-unital \calg $A_0$. Put
$A_{00}=K(A_0)$. Recal that for each $a\in A$, $L_a=cl(A_0 a),
R_a=cl(aA_0)$, and $I_a=cl(L_a R_{a^*})$, where closures are taken
in the norm topology of $A_0$. Also recall that for linear maps
$S:L_a\to L_a $ and $T:R_{a^*}\to R_{a^*}$, $$ M_a=\{(S,T):
yS(x)=T(y)x for each  x\in L_a, y\in R_{a^*}\} $$ is a \calg and
$M_a\cong M(I_a)$, for each $a\in A_0^+$. If $0\leq a\leq b$
are in $A_0$ then $L_a\subseteq L_b, R_a\subseteq R_b,
I_a\subseteq I_b$ and the restriction map defines a \shomo from
$M_b$ to $M_a$. Also for $K_+=K_+(A_0)$ $$ \bigcup_{a\in K_+}
I_a\supseteq \bigcup_{a\in K_+} L_a=\bigcup_{a\in K_+} R_a=K(A_0)\
\  [Ph88b]. $$

Also $K(A_0)$ is a minimal dense ideal, so for each \ai $(e_n)$ of
$A_0$ we have $K(A_{0})\subseteq \bigcup_{n} I_{e_n}$, and so
$$A_0\cong\varinjlim_{n} I_{e_n}, $$ as \calgs. Now let
$A=\Gamma(K(A_{0}))$. Then $A$ is a \sig, and $$
A\cong
\varprojlim_{n} M(I_{e_n}), $$ as \sigs. Also
$A_b=b(A)=M(A_0)$ [Ph88b].

Now, using Dauns-Hofmann's theorem and its generalization to \sigs
, we have that if $X=Prim(A_0)$, then $$ Z(A_b)\cong C_b(X),
~Z(A)\cong C(X) $$ Of course, $Z(A_0)=C_0(X)$ and
$Z(A_{00})=C_{00}(X)$ does not hold in general.

Next let's consider the reverse situation. We want to start with a
\, \,\sig \Asig , with corresponding morphisms $\pi_n : A\to A_n$, and
get a \calg $A_0$ such that $A\cong\Gamma(K(A_0))$. We do
this when $A$ is \lcpt. In this case there is a support algebra
$A_{00}$ in $A$ such that if $A_0$ be the norm closure of $A_{00}$
inside $A_b=b(A)$, then $K(A_0)=A_{00}$. Moreover, for the
corresponding m.b. canonical net $(e_n)$ in $A$, we show that
$(e_n)$ forms an \ai of $A_0$ in norm topology and so by [Ph88b,
thm4] we have $$ \Gamma(A_{00})\cong\varprojlim_{n} M(I_{e_n}),
$$ as \sigs . Finally we show that indeed $$
\Gamma(A_{00})\cong A, $$ as \sigs .

Let's us begin with two well known results about inverse limits.
These are true for uncountable systems (in more general
categories) also, but we state them only for the countable case.
The proofs are quite standard. The first lemma is an immediate
consequence of the definition of inverse limit. We sketch the
proof of the second lemma.
\label{first}
\begin{lemm} If $A=\varprojlim_n (A_n , \pi_n)$ is a \sig and
$\pi_n(A)$ is $C^*$-subalgebra of $A_n$, then
$A\cong\varprojlim_n (\pi_n(A_n) , \pi_n)$, as \sigs .
\end{lemm}

\label{second}
\begin{lemm} If $A=\varprojlim_n (A_n , \phi_n)$ and $B=\varprojlim_n
(B_n , \psi_n)$ are \sigs with connecting morphisms $\phi_{n,n+1}:
A_{n+1}\to A_{n} $ and $\psi_{n,n+1}: B_{n+1}\to B_{n}$, and for
each $n$, there are morphisms $\alpha_n : A_n \to B_n$ and
$\beta_n : B_{n+1} \to A_n$
such that
\[\alpha_n\phi_{n,n+1}=\psi_{n,n+1}\alpha_{n+1} \]
which make the following diagram commute
then $A\cong B$, as topological $*$-algebras.
\end{lemm}
{\bf Proof} By the universal property of inverse limits, there are
maps $j: A\to B$ and $i: B\to A$ such that
for each $n$,
\[\psi_n j=\alpha_n\phi_n\]
and
\[\phi_n i=\alpha_n\beta_n\psi_{n+1}\] and
\[\alpha_n\beta_n=\psi_{n,n+1}\]
the following diagrams commute.
hence $$ \psi_n ji=\alpha_n\phi_n
i=\alpha_n\beta_n\psi_{n+1}=\psi_{n,n+1}\psi_{n+1}=\psi_n , $$ for
each $n$, and so $ji=id_B$. Similarly we get $ij=id_A$ and the
proof is complete.\qed

\vspace{.3 cm}
Since $A_{00}$ is an ideal of $A$, each element of $A$ is a
multiplier of $A_{00}$. This way $A$ maps into $\Gamma(A_{00})$.
The morphism $j:A\to \Gamma(A_{00})$ maps $a$ to the multiplier
$(S_a,T_a)$. Now the compatibility condition tells us that $j$ is
a continuous embedding.

\begin{lemm} If $A$ is \lcpt , then $j$ is one-one.
\end{lemm}
{\bf Proof} If $x\in A$ and $ax=0$ for each $a\in A_{00}$, then in
particular $e_n bx=0$, for each $n$ and $b\in A$. Hence by the
compatibility condition, $\pi_n(x)=0$. Having this for each $n$ we
get $x=0$.\qed

\vspace{.3 cm}
Now we can identify $A$ with a subalgebra of $\Gamma(A_{00})$. Let
the mappings $\pi_n :A\to A_n$ and $\pi_n ^{'}:\Gamma(A_{00})\to
M_{e_n}$ be the corresponding morphisms, where the second is as in
the Theorem 2.2 (denoted by $\pi_n$ there).

\label{weaker}
\begin{prop} With the above notations, if $A$ is locally compact then
$$ \|\pi_{n}(x)\|_{n}\leq\|\pi_n
^{'}(x)\|_n\leq\|\pi_{n+1}(x)\|_{n+1}\ \ (x\in A). $$ In
particular, the embedding $j:A\to\Gamma(A_{00})$ is \cnt and has
continuous left inverse with respect to the corresponding
projective topologies.
\end{prop}
{\bf Proof} Fix $n$ and $x\in A$. Recall from the proof of Theorem
2.2 that $$ \|\pi_n ^{'}(x)\|_n=sup\{\|xae_n\|_\infty
: a\in \Gamma(A_{00}),\ \|ae_n\|_\infty\leq 1\}\geq\|\pi_n(x)\|_n,$$
by the compatibility condition. Let $a\in \Gamma(A_{00})$ be
such that $\|ae_n\|_\infty\leq 1$. Then by the compatibility
condition
\begin{align*}
\|xae_n\|_\infty &=\|xae_ne_{n+1}\|_\infty\\
&\leq\|\pi_{n+1}(xae_n)\|_{n+1}\\
& \leq\|\pi_{n+1}(x)\pi_{n+1}(ae_n)\|_{n+1}\\
& \leq\|\pi_{n+1}(x)\|_{n+1} \|\pi_{n+1}(ae_n)\|_{n+1}\\
& \leq\|\pi_{n+1}(x)\|_{n+1} \|ae_n\|_\infty\\
& \leq\|\pi_{n+1}(x)\|_{n+1}
\end{align*}
which implies that $\|\pi_n ^{'}(x)\|_n\leq
\|\pi_{n+1}(x)\|_{n+1}$ . \qed

\vspace{.3 cm}
It is not easy to show that $j$ is onto. Instead we construct
another morphism from $\Gamma(A_{00})$ to $A$ indirectly. For this
purpose we need an slightly dif and only iferent version of
Phillips' result. First some preliminary results.

Recall that $$ \Gamma(K(A_0))\simeq \varprojlim_{n\in \Bbb N}
(M(I_{e_n}), \pi_{n}^{'})\simeq \varprojlim_{n\in \Bbb N}
(M_{e_n}, \phi_n), $$ where mappings $\pi_{n}^{'}$ are defined by
$\pi_{n}^{''}=\phi_n\circ\psi_n ^{-1}$, where $$
\phi_n:\Gamma(K(A_0))\to M_{e_n}\ \ \text{and}\ \
\psi_n:M(I_{e_n})\to  M_{e_n}, $$ are both defined by restriction
, and the second is an isomorphism [Ph88b]. Now let
$B_n=\phi_n(A_0)$. Then $B_n$ is a $C^*$-subalgebra of $M_{e_n}$
and $\phi_n :A_0\to B_n$ is a surjective \shomo of \calgs , in
particular $\phi_n(K(A_0))=K(B_n)$ [Pd66, II].

The \calgs $B_n$ also provide us another inverse system. Indeed
the restriction of $\phi_{nm}$ to $B_m$ (still denoted by
$\phi_{nm}$) gives us a surjective morphism $\phi_{nm}: B_m\to
B_n$ which extends uniquely to a surjective morphism $\phi_{nm}:
M(B_m)\to M(B_n)$. Also for each $n$, $\phi_n(\Gamma(K(A_0)))=
M(B_n)$ [LT,5.4], hence by Lemma 3.1 we have

\begin{prop} With the above notations,
$$ \Gamma(K(A_0))\simeq \varprojlim_{n} (M(B_n), \phi_{n}), $$ as
\sigs , such that the corresponding morphisms $\phi_n
:\Gamma(K(A_0)) \to M(B_n)$ are all surjective.
\end{prop}

We emphasize that the advantage of this presentation is in the
fact that the morphisms $\phi_n$ are all surjective. Now recall
that $\pi_n ^{'}: \Gamma(K(A_0))\to M(I_{e_n})$ is defined by
$\pi_n ^{'}=\psi_n ^{-1}\circ\phi_n$, where $\psi_n: M(I_{e_n})\to
M_{e_n}$, defined by restriction, is a \calg isomorphism (and so
an isometry). Therefore by Proposition 3.1 we have

\begin{prop} With the above notations, for each $n$ we have
$$ \|\pi_n ^{'}(x)\|_n =\|\phi_n (x)\|_n \ \ \
(x\in\Gamma(K(A_0))), $$ and
 $$ \|\pi_n (a)\|_n
\leq\|\phi_n(a)\|_n\leq\|\pi_{n+1} (a)\|_{n+1} \ \ (a\in A), $$
where in the last equation we have identified $A$ with its image
under $j:A\to \Gamma(K(A_0))$.\qed
\end{prop}

In particular the $\|.\|_\infty$ of $\Gamma(A_{00})$ restricted to
$A$ coincides with the $\|.\|_\infty$ of $A$, namely
 $$ sup_n
\|\pi_n (a)\|_n = sup_n \|\phi_n (a)\|_n\ \ (a\in A). $$

The following lemma is the only place we use the fact that a \lcpt
\sig $A$ is complete with respect to the corresponding
$A_{00}$-topology. Recall that the elements the multiplier algebra
of a Banach algebra $B$ are authomatically bounded and we used the
notation $M(B)$ for the multiplier algebra. This is not the case
when $B$ is merely a normed algebra. Therefore we use the
notations $\Delta(B)$ and $\Gamma(B)$ to distinguish between the
bounded multiplier algebra and the multiplier algebra.

\label{didit}
\begin{lemm} With the above notations, if $A$ is \lcpt and
$A_{00}$ is the corresponding support algebra of $A$, then
$\Delta(A_{00})\subseteq A\subseteq \Gamma(A_{00})$.
\end{lemm}
{\bf Proof} Let $x\in \Delta(A_{00})$ be self adjoint and let
$(e_n)$ be as in the Definition 2.6. Then we claim that the
sequence $(e_n x)\subseteq A$ is $A_{00}$-Cauchy. Note that, since
$x$ is a multiplier of $A_{00}$ and $e_n\in A_{00}$, we have $e_n
x\in A_{00}\subseteq A$. Now, for each $a\in A_{00}$, we have $$
\|a(e_n x-e_m x)\|_\infty\leq\|a(e_n -e_m
)\|_\infty\|x\|_\infty\to 0 , $$ and
 $$ \|(e_n x-e_m
x)a\|_\infty\leq\|e_n(xa) -e_m(xa)\|_\infty\to 0 , $$ as
$m,n\to\infty$. All of this calculation is done inside
$\Gamma(A_{00})$. In particular $\|.\|_\infty =sup_n \|\phi_n
(.)\|_n$. But from the above lemma this calculation is valid
inside $A$ also (because $ sup_n \|\pi_n (.)\|_n = sup_n \|\phi_n
(.)\|_n$). This completes the proof of the claim.

Now, as $A$ is $A_{00}$-complete, there is a self adjoint element
$y\in A$ such that $e_n x\to y$ in $A$ with respect to the
$A_{00}$-topology. In particular, for each $a\in A_{00}$, we have
$\|e_n xa-ya\|_\infty\to 0$ in $A$. On the other hand, $\|e_n
xa-xa\|_\infty \to 0$ in $\Gamma(A_{00})$. Again all these nets
are inside $A$ and the $\|.\|_\infty$ of $\Gamma(A_{00})$
restricted to $A$ coincides with $\|.\|_\infty$ of $A$. Therefore
$ya=xa\in A_{00}$. But $x$ and $y$ are self adjoint and $A_{00}$
is an $*$-algebra, hence $ay=ax\in A_{00}$. Therefore $x=y$ as
elements of $\Gamma(A_{00})$ and the proof is finished.\qed

\vspace{.3 cm}
The first part of the following lemma is a modification of [LT,
5.4].

\label{phi}
\begin{lemm} With the above notations, for each n
$$ \phi_n(\Delta(K(A_0)))=\phi_n(\Gamma(K(A_0)))=\Gamma(K(B_n))=
\Delta(K(B_n))\subseteq M_{e_n}. $$ Moreover if $A$ is \lcpt ,
then we have $j(A)\subseteq \Gamma(K(A_0))$ and
$\phi_n(j(A))=\Delta(K(B_n))=M(B_n)$.
\end{lemm}
{\bf Proof} Recall that $A_0$ (and so its quotient $B_n$) is a
$\sigma$-unital \calg . Now the second equality and the inclusion
holds for any \calg $A_0$ by [LT, 5.4]. The first equality is
proved in [LT, 5.4] under the assumption that $A_0$ is separable.
However, they use this assumption only to make sure that the
surjection $\phi_n :A_0\to B_n$ extends to a surjective morphism
between the multiplier algebras ( see [APT, 4.2] quoted in [LT,
5.4] as Theorem 9.2 which is a typo!). But this holds also for
$\sigma$-unital \calgs (this is a special case of [Ph88a, 5.11]).

Now let $A$ be \lcpt and $A_{00}$ be the support algebra of
$A$(which we have chosen and fixed), then we know that
$j(A)\subseteq \Gamma(A_{00})$. Hence $\phi_n(j(A))\subseteq
\Delta(K(B_n))$ by the first part. Now by the above paragraph,
$\phi_n(M(A_0))=M(B_n)$, hence by Lemma 3.4, we have
$\phi_n(j(A))\supseteq \phi_n(\Delta(K(A_0)))=\Delta(K(B_n))$.
Also $\phi_n(j(A))\subseteq
\phi_n(\Gamma(K(A_0)))=\Delta(K(B_n))$, and we are done.\qed

\vspace{.3 cm}
Now recall that if $A$ is \lcpt then the morphisms $\pi_n :A\to
A_n$ are surjective (c.f. Definition 2.6). Also in this case Lemma
3.5 tells us that $\phi_n\circ j :A\to M(B_n)$ is surjective for
each $n$. Now by the same lemma the morphisms
$$
\tilde\pi_n:A_{n+1}\to M(B_n)\ \ \ \pi_{n+1}(a)\mapsto \phi_n j(a)
$$
and
$$
\hat\pi_n : M(B_n)\to A_{n}\ \ \ \phi_n j
(a)\mapsto\pi_{n}(a)
$$
are well defined, norm \cnt , and
surjective. Now we are ready to prove the main result of this
section.

\label{main}
\begin{theo}  Let $A=\varprojlim_{n\in \Bbb N} A_n$ be a locally compact
\sig . Let $(e_n)\subseteq A_{sb}$ and
$A_{00}=\bigcup Ae_n ^2 A$ be the corresponding approximate
identity and support algebra. Let $A_0$ be the C$^*$-subalgebra of
$A_b$ which is the norm closure of $A_{00}$ in $A_b$. Then
$A_{00}=K(A_0)$ and $A\cong\Gamma(K(A_{0}))$, as \sigs , and
$A_b\cong M(A_0)$, as \calgs. Moreover the \calg $A_0$ is
unique if we require that the isomorphism from $A$ onto
$\Gamma(K(A_0))$ is $A_{00}$-$\kappa$-bicontinuous.
\end{theo}
{\bf Proof} The fact that $A_{00}=K(A_0)$ is Lemma 2.5. Now from
the two last isomorphisms, the second follows immediately from the
first and the fact that $\Gamma(K(A_{0}))_b =M(A_0)$ [Ph88b]. For
the first statement, by Lemma 3.2, we only need to observe that,
for each $n$,
$$
\pi_{n,n+1} \tilde\pi_{n-1}=\phi_{n-1,n}\hat\pi_n
$$
But this follows directly from the definition of the maps
$\tilde\pi_{n}$ and $\hat\pi_{n}$, and that each map $\phi_n j$
are onto (by Lemma 3.5). For the last statement, let $A_0$ and
$B_0$ be two \calgs such that $\Gamma(K(A_0))$ and
$\Gamma(K(B_0))$ are both isomorphic to $A$ with the isomorphisms
being $A_{00}$-$\kappa$ and $B_{00}$-$\kappa$ bicontinuous,
respectively. Then the composition of these would be a
$\kappa$-$\kappa$-bicontinuous isomorphism between
$\Gamma(K(A_0))$ and $\Gamma(K(B_0))$. Then $A_0$ and $B_0$ are
isomorphic by [LT, thm 7.10].\qed
\begin{rem} It is clear from the above proof that the {\it
compatibility condition} in the definition of \lcpt \sigs could be
replaced with the following
$$
\|e_{n_{k-1}}a\|_\infty\leq\|\pi_k(a)\|_k\leq
sup\{\|e_{n_k}ba\|_\infty:\ b\in A, \|e_{n_k}b\|_\infty \leq 1\}\
\ \ (a\in A, k\in\Bbb N),
$$
 for some (infinite) subsequence
$\{n_k\}$ of $\Bbb N$. This, in particular, justifies Example 2.1.
\end{rem}

The above theorem with the previously mentioned result of N.C.
Phillips give a correspondence between the {\it objects} of the
given categories. Namely

\begin{cor} A \sig $A$ is \lcpt if and only if there is a \calg $A_0$ such that
$A \cong\Gamma(K(A_{0}))$, as \sigs .
\end{cor}
{\bf Proof} If $A$ is \lcpt and $A_0$ is as in the above theorem,
then $A\cong \Gamma(K(A_{0}))$. The other direction is
Theorem 2.2.\qed

\vspace{.3 cm}
If $A=C(X)$ is a \lcpt commutative \sig , where $X$ is a \hau
$\kappa_\omega$-space, then the above result tells us that there
exists a $C^*$-subalgebra $A_0$ of $A_b$ such that $A
\cong\Gamma(K(A_{0}))$, as \sigs . But $A_0=C(Y)$, for a
\lcpt $\sigma$-\cpt \hau space $Y$ and so $C(X)$ and $C(Y)$ are
isomorphic. In general this does {\it not} imply that $X$ and $Y$
are homeomorphic. But here both $X$ and $Y$ are $\sigma$-\cpt and
so they are real compact (see [GJ, 8.2] for the definition and
proof). In particular, $X$ and $Y$ are homeomorphic [GJ, 10.6].
Hence $X$ is \lcpt . This gives an alternative proof of (the
difficult direction of) Theorem 2.1. Also it shows that the
uniqueness part of the Theorem 3.1 holds in the commutative case.

\begin{cor} If \Asig is a \lcpt \sig and $A_{00}$ is
support algebra of $A$, then the projective topology is stronger
than the $A_{00}$-topology.\qed
\end{cor}

Next let us consider the issue of morphisms. Let $A_0, B_0$ be
\calgs and $\phi_0: A_0\to B_0$ be a \shomo, then $\phi_0$ does
not lift to a \shomo of the corresponding multiplier algebras
unless it is surjective. To turn around this difficulty people
usually consider
 $$ Mor(A_0 ,B_0)=\{\phi_0 :A_0\to M(B_0):\
\phi_0(A_0)B_0\subseteq B_0 \ \ \text{is \ dense}\} $$ as the
family of {\it morphisms} from $A_0$ to $B_0$. Such morphisms are
called {\it non degenerate}. These morphisms have a unique
extension to (strictly continuous) unital $*$-homomorphisms
between multiplier algebras.

\begin{defi} Let $A_0$ and $B_0$ be as above, a morphism $\phi_0
:A_0\to M(B_0)$ is called strictly non degenerate if
$\phi_0(A_{00})B_0\supseteq B_{00}$, where $A_{00}$ and $B_{00}$
are the Pedersen ideals of $A_0$ and $B_0$, respectively.
\end{defi}

Now let $A, B$ be the corresponding unital \sigs, i.e.
$A=\Gamma(A_{00})$ and $B=\Gamma(B_{00})$, where $A_{00}$ and
$B_{00}$ are the Pedersen ideals of $A$ and $B$, respectively. A
{\it morphism} from $A$ to $B$ is a unital \shomo such that
 $$
\forall m\ \exists n\ \ \ \|\phi(a))\|_m\leq\|a\|_n, \ \ (a\in A).
$$
This condition plays a crutial role. It implies that each \repn of a \sig \Apro factors
through some $A_n$. It also ensures that the GNS-construction
yields a \repn of $A$.

\label{ndeg}
\begin{defi}
We say that $\phi$ is {\it non degenerate} if
$\phi(A_{00})B\subseteq B_0$ is dense. We say that $\phi$ is {\it
strictly non degenerate} if $\phi(A_{00})B\supseteq B_{00}$.
\end{defi}

Now let us start with a morphism $\phi_0\in Mor(A_0,B_0)$. We wish
to extend $\phi_0$ to a morphism $\phi:A\to B$.

\label{extend}
\begin{lemm} 
With the above notation, if $\phi_0$ is (strictly) non degenerate
then we have $B_0\phi_0(A_{00})B_0=B_{00}$
\,\,(\, $B_0\phi_0(A_{00})=\phi_0(A_{00})B_0=B_{00}$, respectively).
\end{lemm}
{\bf Proof} Since $\phi_0$ preserves the spectral theory,
$\phi_0(A_{00})\subseteq B_{00}$. On the other hand,
$\phi_0(A_{00})$ is clearly a dense ideal of $\phi_0(A_0)$.
Therefore $B_0\phi_0(A_{00})B_0$ is an ideal of $B_0$ which is
dense in $B_0\phi_0(A_0)B_0$. Now if $\phi_0$ is non degenerate
then $\phi_0(A_0)B_0$ is dense in $B_0$, so $B_0\phi_0(A_0)B_0$ is
dense in $B_0 B_0 =B_0$, i.e. $B_0\phi_0(A_{00})B_0$ is a dense
ideal of $B_0$ and so contains $B_{00}$. But
$B_0\phi_0(A_{00})B_0\subseteq B_0 B_{00} B_0\subseteq B_{00}$,
hence the equality holds. If $\phi_0$ is strictly non degenerate
then $\phi_0(A_{00})B_0\supseteq B_{00}$. The converse inequality
follows from the fact that $\phi_0(A_{00})\subseteq B_{00}$. Hence
$\phi_0(A_{00})B_0=B_{00}$. Now the right hand side is self
adjoint and the adjoint of the left hand side is
$B_0\phi_0(A_{00})$, hence $B_0\phi_0(A_{00})=B_{00}$.\qed

\begin{rem}
When $\phi_0$ is surjective, one  gets the better result
$\phi_0(A_{00})=B_{00}$, first proved by G.K. Pedersen.\
\end{rem}

Now let $\phi_0$ be strictly non degenerate. Then it extends to a
map $\phi:A=\Gamma(A_{00})\to B=\Gamma(B_{00})$ given by
$$
 \phi(x)\phi_0(a)b=\phi_0(xa)b,\,\, b\phi_0(a)\phi(x)=b\phi_0(ax)
 \quad (x\in A, a\in A_{00}, b\in B_0).
$$
Again, as $\phi$ preserves the spectral theory, we have
$\phi(A_b)\subseteq B_b$, and so $\phi$ induces a $*$-homomorphism
$\phi_b:A_b\to B_b$. By uniqueness, $\phi_b$ is the same as the
extension of $\phi_0$ to $A_b=M(A_0)$. Next observe that $\phi$ is
strictly non degenerate. Indeed
$\phi(A_{00})B\supseteq\phi_0(A_{00})B_0\supseteq B_{00}$.
Now we want to examine the other direction. This time a strictly
non degenerate morphism $\phi:A\to B$ is given, and we are aiming
to show that there exist a strictly non degenerate morphism
$\phi_0:A_0\to B_0$ which coincides with the restriction of $\phi$
to $A_0$. Indeed the fact that $\phi$ is a morphism in particular
implies that $\phi(A_b)\subseteq B_b$, so if we take $\phi_0$ to
be the restriction of $\phi$ to $A_0$, then $\phi_0$ is norm
continuous (indeed of norm$\leq 1$) and $\phi_0(A_0)\subseteq
\phi(A_b)\subseteq B_b$.

Now $\phi$ preserves the spectral theory, hence
$\phi(A_{00})\subseteq B_{00}$, and so $\phi(A_{00})B\subseteq
B_{00}B\subseteq B_{00}$. Hence
$\phi(A_{00})B=\phi(A_{00})B_0=B_{00}$. We need to show that
$\phi_0(A_{00})B_0\supseteq B_{00}$. But since $\phi_0$ is the
restriction of $\phi$ to $A_0$, we have $\phi_0(A_{00})B_0
=\phi(A_{00})B_{0}=B_{00}$. Therefore we have shown that, given a
strictly non degenerate morphism $\phi:A\to B$, of \sigs the
restriction $\phi_0$ of $\phi$ to $A_0$ is a strictly non
degenerate morphism of \calgs . Conversely, each strictly non
degenerate element $\phi_0\in Mor(A_0,B_0)$ uniquely extends to a
strictly non degenerate morphism $\phi:A\to B$.


Let's summarize these observations as follows.

\label{cateq}
\begin{theo}
There is a covariant functor from the category of $\sigma$-unital
\calgs and strictly non degenerate $*$-homomorphisms to the
category of \lcpt \sigs and strictly non degenerate
$*$-homomorphisms. This functor assigns to each $\sigma$-unital
\calg the \sig of multipliers of its Pedersen ideal. In
particular, for each $\sigma$-\cpt \hau \ts $X$, it sends the
commutative \calg $C_0(X)$ to the commutative \sig $C(X)$.
\end{theo}
{\bf Proof} We have established a one to one correspondence
between the objects and morphisms of the category of locally
compact \sigs and the category of \calgs. The fact that this is
indeed an equivalence of categories then follows from the trivial
observation that the following diagram commutes:
\[\iota_B\phi_0=\phi\iota_A\]
The other statements are trivial.\qed

\vspace{.5 cm}
\section{Examples}

In this section we use the results of previous sections to
calculate the \pros associated with some of the well known \calgs .

\begin{ex} Let $X$ be a \lcpt \hau \ts, and $A_0=C_0(X)$, then
$A_{00}=K(A_0)=C_{00}(X)$, and so $A=\Gamma(A_{00})=C(X)$, and
$A_b=M(A_0)=C_b(X)$. Conversely, if $A=C(X)$ is \lcpt as a \pro,
then $X$ is \lcpt as a \ts  and $A_b=C_b(X)$. Also $A_{00}$
consists of \cnt \fns whose multiplication with any \cnt \fn is
bounded. Such functions are exactly the ones of compact support.
Hence $A_{00}=C_{00}(X)$ and therefore $A_0=C_0(X)$.
\end{ex}

\begin{ex} Let \Apro and \Bpro be \pros. Recall that $$
A\bigotimes_{max} B=\varprojlim_{(\alpha,\beta)}
A_\alpha\bigotimes_{max} B_\beta, ~~A\bigotimes_{min}
B=\varprojlim_{(\alpha,\beta)} A_\alpha\bigotimes_{min} B_\beta.
$$ Let $A_0$ and $B_0$ be the corresponding \calgs, then \newline
\underbar{Claim}: $(A\bigotimes_{max} B)_0=A_0\bigotimes_{max}
B_0$ and $(A\bigotimes_{min} B)_0= A_0\bigotimes_{min}
B_0$.\newline Let $(e_\alpha)\subseteq A_{00}$ and
$(f_\beta)\subseteq B_{00}$ be approximate identities of $A$ and
$B$, respectively. Then $(e_\alpha\otimes f_\beta)\subseteq
A_{00}\bigodot B_{00}\subseteq (A\bigotimes B)_{00}$, where
$\bigotimes$ is any of the $max$ or $min$ tensor products. Let
$I_\alpha=I_{e_\alpha}$ and $J_\beta=I_{f_\beta}$. Then $$
I_\alpha\bigotimes J_\beta\cong I_{(e_\alpha\otimes f_\beta)}
$$ as \calgs. Indeed the left hand side is included in $ A_0e_\alpha A_0\bigodot B_0
f_\beta B_0=(A_0\bigodot B_0)(e_\alpha\otimes f_\beta)(A_0\bigodot
B_0)$, which is norm dense in the right hand side. Also the right hand side is included in
$(A\bigotimes B)_0(e_\alpha\otimes f_\beta)(A\bigotimes
B)_0\supseteq (A_0\bigodot B_0)(e_\alpha\otimes
f_\beta)(A_0\bigodot B_0)= A_0e_\alpha A_0\bigodot B_0 f_\beta
B_0$, which is norm dense in the left hand side.  Now direct product
preserves the $C^*$-tensor products, so 
\begin{align*}
 A_0\bigotimes
B_0 &\cong \varinjlim_{\alpha}I_\alpha\bigotimes
\varinjlim_{\beta} J_\beta\cong \varinjlim_{(\alpha,\beta)}
I_\alpha\bigotimes J_\beta\\
&\cong \varinjlim_{(\alpha,\beta)}
I_{e_\alpha\otimes f_\beta}\cong (A\bigotimes B)_0. 
\end{align*}
Note that the similar result for $(A\bigotimes B)_{00}$ needs an
appropriate topological tensor product completion of
$A_{00}\bigodot B_{00}$. Even in special cases this seems to be
unavailable (see example ... below). Also even in commutative
case, $(A\bigotimes B)_b\varsupsetneqq A_b\bigotimes B_b$.
\end{ex}

\begin{ex} Let $A_0=K(H)$ be the algebra of all \cpt operators on
a Hilbert space $H$. Then $A_{00}=F(H)$ is the algebra of all
operators of finite rank. Hence $A=\Gamma(F(H))=B(H)$. In this
case $A_b=A=B(H)$. Also if $A_0=B(H)$ then
$A_{00}=A_0=A_b=A=B(H)$.
\end{ex}

\begin{ex} Let $X$ be a \lcpt \hau \ts, and $A_0$ be a (non
unital) \calg. Then $B_0=C_0(X,A_0)\cong
C_0(X)\bigotimes_{min} A_0$ is also a \calg whose multiplier
algebra is $B_b=M(C_0(X,A_0))=C_b(X,A_b)$ [Wr95]. But $B_{00}\ne
C_{00}(X,A_{00})$. Indeed if $A_0=K(H)$ then $A_{00}=F(H)$ is the
algebra of all operators of finite rank, and $A=B(H)$. Now
$B=C(X)\bigotimes_{min} A$ is the tensor product of two \lcpt
\pros, and so is \lcpt. Therefore $B_{00}=K(B_0)=K(C_0(X,K(H)))$
which is much smaller than $C_{00}(X,F(X))$ (it is even a proper
subset of the set of all those \cnt \fns $f:X\to F(X)$ for which
$\sup_{x}(dim f(X))<\infty$) [GT].
\end{ex}

\begin{ex} Let $G$ be a \lcpt \hau topological group acting on a
\calg $A_0$ via a \cnt group \homo $\alpha_0:G\to Aut(A_0)$. Let
$A_{00}=K(A_0)$, $A_b=M(A_0)$ and $A=\Gamma(A_{00})$. For each
$g\in G$, the corresponding $*$-automorphism $\alpha_{0}(g):A_0\to
A_0$ is in particular surjective. Therefore both of its
restriction $\alpha_{00}(g): A_{00} \to A_{00}$, and its extension
$\alpha_b(g):A_b\to A_b$ are also surjective. Indeed they are also
injective. This is trivial for the first map. For the second,
let's recall that each $x\in A_b$ could be considered as an
element of $A_{0}^{**}$ such that $xA_0\cup A_0x\subseteq A_0$.
Let's fix $g\in G$. We abbreviate $\alpha_0(g)(a)=g.a, a\in A_0$.
Then for each $x\in A_b$, define $g.x\in A_b=M(A_0)$ by
$(g.x)a=g.(x(g^{-1}.a)), a\in A_0$. This agrees with the our
previous notation when $x\in A_0$, therefore, by uniqueness of
extension, $\alpha_b(g)(x)=g.x$, which is clearly injective.

Similarly $\alpha_0(g)$ extends to a to an automorphism of $A$
given by $$(g.x)a=g.(x(g^{-1}.a)), \quad (a\in A_{00}, x\in A).$$ This
defines an action $\alpha:G\to Aut(A)$. Now if $A=\varprojlim_{i}
A_i$, where $A_i$'s are (unital) \calgs and {\it all morphisms
$\pi_i: A\to A_i$ are surjective}, then $G$ acts on $A_i$ via
$g.\pi_{i}(x)=\pi_{i}(g.x); ~g\in G, x\in A$. Let's observe that
the actions $\alpha_i$ of $G$ on $A_i$ are compatible with the
inverse system, i.e. $\alpha_i(g)\pi_i=\pi_i\alpha_i(g)$,  for
each $i$, and each $g\in G$, which is just the definition of
$\alpha_i$ (note that $\alpha_i$ is well defined, because
$\pi_i$'s are surjective). Therefore it is reasonable to define
the {\it crossed product} of $G$ with $A$ by $$ G\times_{\alpha}
A=\varprojlim_{i} G\times_{\alpha_i} A_i. $$ It would be desirable
to show that $$ G\times_{\alpha} A=\Gamma(K(G\times_{\alpha_0}
A_0)), $$ but this is not true in general (take $A=\Bbb C$).
\end{ex}

\begin{ex} Let $G$ be as above and $A_0=C^*(G)$ be the group
\calg. As far as I know, there is no specific way to describe the
Pedersen ideal of $C^*(G)$ in general. In particular, I don't know
how to calculate the corresponding \pro. However, if $G$ is a {\it
[SIN]-group} (i.e. it has a local basis of neighborhoods of
identity which are invariant under inner automorphisms), then one
has an explicit description of it: Let $\hat{G}$ be the set of
(representatives of the equivalence classes of ) all irreducible
\repns of $G$. Each $a\in C^*(G)$ induces a \cnt \fn
$\hat{a}:\hat{G}\to B(H)$, defined by $\hat{a}(\pi)=\pi(a)$, where
$H$ is the Hilbert space of the universal \repn of $G$. A subset
$K$ of $\hat{G}$ is called {\it quasi-compact}, if each $\hat{a}$
is bounded on $K$. Then the Pedersen ideal of $C^*(G)$ could de
described as $$ K(C^*(G))=\{a\in C^*(G): \hat{a} ~vanishes
~outside ~a ~quasi-compact ~subset ~of ~\hat{G}\}. $$ If we equip
$\hat{G}$ with the {\it quasi-topology}, then the Pedersen ideal
is exactly the ideal which corresponds to the ideal of compactly
supported \fns on $\hat{G}$, i.e. $C_{00}(\hat{G})=K(C^*(G))$, as
sets (after trivial identification). Indeed for each $a\in
C^*(G)$, $\hat{a}$ vanishes at infinity on $\hat{G}$: Given
$\epsilon> 0$, put
$C=\{\pi\in\hat{G}:\|\hat{a}(\pi)\|\geq\epsilon\}$, then $C$ is
quasi-\cpt [Kan] and $\|\hat{a}(\pi)\|< \epsilon$ off $C$. If
$C^*(G)$ is separable, then the converse is also true and
$C^*(G)=C_0(\hat{G})$, as sets [Kan].

Now, take any $a\in K(C^*(G))$, and consider the quasi-\cpt subset
$K_a=\{\pi\in\hat{G}:\|\hat{a}(\pi)\|\geq 1\}$, then
$C_{a}^*(G)=\{b\in C^*(G): \hat{b}\ is \ zero \ off \ K_a\}$ is a
unital \calg, and $A_{00}=K(C^*(G))=\bigcup_{a} C_{a}^*(G)$, where
$a$ runs over $K(C^*(G))$. Therefore $$ A_0=C^*(G)=\varinjlim_{a}
C_{a}^*(G), $$ and $$ A=\Gamma(K(C^*(G)))=\varprojlim_{a}
M(C_{a}^*(G)), $$ where $\pi_{ab}:C_{a}^*(G)\to C_{b}^*(G)$ is
just the inclusion map, when $a\leq b$.
\end{ex}

\begin{ex} Let $(A_\alpha)_{\alpha\in\Lambda}$ be a net of unital
\calgs, and $A_0=\sum_{\alpha}^{c_0} \bigoplus
A_\alpha=\{(a_\alpha)\in \prod_{\alpha} A_\alpha:\lim_{\alpha}
\|a_{\alpha}\|=0\}$ with the norm
$\|(a_{\alpha})\|_{\infty}=\sup_{\alpha}\|a_{\alpha}\|$.
 Then $A_0$ is a (non unital) \calg with the Pedersen ideal
$A_{00}=K(A_0)=\bigcup_{\alpha} (\sum_{\beta\leq\alpha}^{\oplus}
A_{\beta})$. In particular $A_0=\varinjlim_{\alpha}
(\sum_{\beta\leq\alpha}^{\oplus} A_{\beta}).$ Also
$A_b=M(A_0)=\{(a_\alpha)\in\prod_{\alpha} A_\alpha: \sup_{\alpha}
\|a_{\alpha}\|<\infty\}=\sum_{\alpha}^{ l^{\infty}}\bigoplus
A_\alpha$, and $$ A=\Gamma(A_{00})=\prod_{\alpha}
A_\alpha=\varprojlim_{\alpha} (\sum_{\beta\leq\alpha} \bigoplus
A_{\beta}). $$ To prove the above assertions, let's first observe
that, for each $\alpha$,
$B_{\alpha}=\sum_{\beta\leq\alpha}^{\oplus} A_{\beta}=\{(a_\beta):
a_\beta=0, for \beta >\alpha\}$ is a unital $C^*$-subalgebra of
$A_0$. Put $A_{00}=\bigcup_{\alpha} B_{\alpha}$. This is a union
of an increasing family of ideals of $A_0$, so it is an ideal.
Moreover it is norm dense: Given $a=(a_{\alpha})\in A_0$, and
$\epsilon\geq 0$, there is $\alpha_0$ such that
$\|a_{\alpha}\|\leq\epsilon$, ~$\alpha\geq\alpha_0$ . Let
$a^{'}=(a_{\alpha}^{'})$, where $a_{\alpha}^{'}=a_{\alpha}$ for
$\alpha\leq\alpha_0$, and zero otherwise, then
$\|a^{'}-a\|_{\infty}\leq\epsilon$, as required. In particular,
$A_{00}\supseteq K(A_0)$. Conversely, consider the projections
$\pi_\alpha: \prod_{\alpha} A_\alpha\to B_\alpha$ given by
$\pi_\alpha((a_\beta))=(a^{'}_\beta)$, where
$a_{\beta}^{'}=a_\beta$, for $\beta\leq\alpha$, and zero
otherwise, then it is obvious that this is a surjection (even if
it is restricted to $A_0$), so it sends $K(A_0)$ onto
$K(B_\alpha)=B_\alpha$. But this just means $K(A_0)\supseteq
B_\alpha$, for each $\alpha$. In particular, $K(A_0)\supseteq
A_{00}$, so the equality holds.

Next Let $A_b=\sum_{\alpha}^{l^{\infty}} \bigoplus A_\alpha$, then
let's observe that $A_b=M(A_0)$: $A_0$ is clearly an ideal in
$A_b$. It is also essential. Indeed if $a=(a_\alpha)\in A_b$ and
$aA_0=\{0\}$, then, given $\alpha$, there is $b=(b_\beta)\in
A_{00}$ such that $b_\alpha=1_\alpha$(=the unit element of
$A_\alpha$), in particular, $ab=0$ implies that $a_\alpha=0$, so
$a=0$. Therefore, there is a canonical embedding $A_b\subseteq
M(A_0)$. \ Conversely, $A_{0}^{**}=(\sum_{\alpha}^{c_0} \bigoplus
A_\alpha)^{**} =\sum_{\alpha}^{l^{\infty}} \bigoplus
A_{\alpha}^{**}$, and each $b\in M(A_0)$ is of the form
$b=(x_\alpha)$, where $x_\alpha\in A_\alpha^{**}$, such that
$bA_0\cup A_0b\subseteq A_0$. In particular, $x_\alpha
A_\alpha\cup A_\alpha x_\alpha\subseteq A_\alpha$, that is
$x_\alpha\in M(A_\alpha)=A_\alpha$. This means that
$M(A_0)\subseteq \sum_{\alpha}^{l^{\infty}} \bigoplus
A_\alpha=A_b$.\newline Finally, $$
A=\Gamma(A_{00})=\Gamma(\bigcup_{\alpha}
B_\alpha)=\varprojlim_{\alpha} (\sum_{\beta\leq\alpha} \bigoplus
A_{\beta})=\prod_{\alpha} A_\alpha. $$ A special case of this
example is when the original $A_\alpha$'s are (a countable family
of) matrix algebras, then $A_0$ is an algebra of {\it block
matrices} of infinite size, $A_{00}$ are those which have only
finitely many nonzero blocks, and $A$ consists of all infinite
matrices. Another special case of this is already considered by
S.L. Woronowicz [Wr91].
\end{ex}

\begin{ex}  Let $p$ and $q$ be the {\it momentum} and {\it
position} operators of a quantum mechanic system of one degree of
freedom. In {\it Schr\"odinger} \repn, $H=L^{2}({\Bbb R},dx)$, and
$p=M_x$, $q=-\frac{d}{dx}$. Then $p,q\eta K(H)$, but $p,q\notin
B(H)=\Gamma(F(H))$.
\end{ex}

{\small {\bf Acknowledgement}: This paper is part of the author's
Ph.D. thesis in the University of Illinois at Urbana-Champaign
under the supervision of Professor Zhong-Jin Ruan. I would like to
thank him for his moral support and scientific guidance during my
studies.}


\begin{thebibliography}{99}      

\bibitem[APT]]
 C.A. Akemann, G.K. Pedersen, J. Tomiyama: Multipliers of \calgs.
J. Functional Analysis 13(3)  (1973) , 277-301.

\bibitem[CD]]
 H.S. Collins, J.R. Dorroh: Remarks on certain function spaces .
Math. Ann. 176  (1968) , 157-168.

\bibitem[CF]]
 H.S. Collins, R.A. Fontenot: Approximate identities and the strict
topology . Pacific j. Math. 43   (1972) , 63-80.

\bibitem[FW]]
 R.A. Fontenot, R.F. Wheeler: Approximate identities and
paracompactness. Proc. Amer. Math. Soc. 99  (1987) , 232-236.

\bibitem[GJ]]
 L. Gillman, M. Jerison: Rings of \cnt \fns (1976) , Springer
Verlag , New York, Heidelberg, Berlin.

\bibitem[Lin]]
 H.X. Lin: Support algebras of $\sigma$-unita \calgs and
their quasi-multipliers . Trans. Amer. Math. Soc.
 325  (1991) , 829-854.

\bibitem[LT]]
 A. Lazar, D.C. Taylor: Multipliers of Pedersen's ideal . Mem.
Amer. Math. Soc. 169  (1976).

\bibitem[Pd66]]
 Gert K. Pedersen: Measure theory for \calgs $I-IV$. Math. Scand.
19  (1966) , 131-145; {\bf 22} (1968) 63-74; {\bf 25} (1969)
71-93; {\bf 25} (1969) 121-127.

\bibitem[Pd64]]
 Gert K. Pedersen: Measure theory for \calgs , Ph.D. Thesis.

\bibitem[Ph88a]]
 N.C. Phillips: Inverse limits of \calgs. J. Operator Theory
 19  (1988) , 159-195.

\bibitem[Ph88b]]
 N.C. Phillips: A new approach to the multipliers of Pedersen's
ideal . Proc. Amer. Math. Soc. 104(3)  (1988) , 861-867.

\bibitem[Ph88c]]
 N.C. Phillips: Inverse limits of $C\sp *$-algebras and
applications, Operator algebras and applications, Vol. 1, 127-185,
London Math. Soc. Lecture Note Ser., 135 (1988) , Cambridge Univ.
Press , Cambridge-New York.

\bibitem[Ty72]]
 D.C. Taylor: A general Phillips theorem for \calgs and some
applications . Pacific J. Math. 40  (1972) , 477-488.

\bibitem[Whe]]
 R.F. Wheeler: Well-behaved and totally bounded approximate
identities for $C_0(X)$ . Pacific J. Math. 65  (1976)  , 261-269.

\end{thebibliography}
\end{document}